# New Geometric Constant Related to the *P*-angle Function in Banach Spaces

Zhijian Yang and Yongjin Li *

Department of Mathematics, Sun Yat-sen University, Guangzhou, 510275, P. R. China

**Abstract**

In this paper, combined with the *P*-angle function of Banach spaces and the geometric constants that can characterize Hilbert spaces, the new angular geometric constant $S_P(X)$ is defined. Firstly, this paper explores the basic properties of the new constant $S_P(X)$ and obtains some inequalities with significant geometric constants. Then according to the derived inequalities, this paper studies the relationship between $S_P(X)$ and the geometric properties of Banach spaces. Furthermore, the necessary and sufficient condition for uniform non-squareness, and the sufficient conditions for uniform convexity, the normal structure and the fixed point property will be established.

**Mathematics Subject Classification:** 46B20

**Keywords:** Banach spaces; P-angle function; geometric constants; the fixed point property

# 1 Introduction

As we all know, the geometric theory of Banach spaces has been fully developed and synthesizes the properties of concrete spaces such as the classical sequence spaces $c_0$, $l_p(1 \leq p < \infty)$ and the function space $C[a, b]$. In particular, the geometric constant has become a powerful tool to characterize the geometric properties of the space sphere. After fifty years of exploration and research, scholars found that some abstract properties of Banach spaces can be quantitatively described by some special constants. At present, there are many papers and achievements on geometric constants, but how to use the geometric constants to classify Banach spaces is an important problem. For example, Clarkson introduced the module of convexity to be used to characterize uniformly convex spaces [15], and the von-Neumann constant to be used to characterize uniformly non-square spaces [6]. In order to study the normal structure of spaces, James introduced the James constant [13]. After the appearance of these constants, many scholars paid attention to them and obtained many wonderful properties. Although the study of geometric constants has gone through more than half a century, many new geometric constants constantly appear in our field of vision. Since the 1960s, not only the geometric theory of Banach spaces has been fully developed, but also its research methods have been applied to matrix theory, differential equations and so on.

In the traditional European spaces, the right angle which refers to the angle in the traditional sense and the orthogonality are closely related. The concept of angle and how to measure angle are ancient and interesting mathematical problems. Obviously, the construction of pyramids and other ancient monuments are inseparable from angle. In the eighteenth and nineteenth centuries, the Euclidean geometry gradually lost its effectiveness. Obviously, in various cases, the concept of angle and the measurement of angle should be reconsidered. Spherical geometry is the first non-Euclidean geometry, in which the internal concept of angle seriously requires us to define and study the angle of spatial curve that is no longer Euclidean line. In 2017, in order to better characterize the properties of inner product spaces, Vitor Balestro [5] defined some new angles, and deduce the theorems of Banach spaces by these angles. In the traditional Euclidean plane, the angle has strong geometric intuition and has been widely accepted, but there is no good conclusion on the characterization of general Banach spaces. Geometric constants can describe general Banach spaces, but also lose some geometric properties. Therefore, in recent years, the idea of exploring the properties of Banach spaces through the geometric constants defined by new angles has been practiced by many scholars, especially the constants defined by the right angles, including *D(X), BR(X), IB(X), D₀(X)* and so on [17].

Although Vitor Balestro [5] introduced some new angles functions, including *P−* angle, *S−*angle and *I−*angle, the work of using these angles to study the geometric properties of Banach spaces has not yielded good results. Minkowski Geometry, the geometry of real finite dimensional Banach spaces, is a fascinating mathematical discipline, which is closely related to other mathematical fields, such as Functional Analysis, Distance Geometry, Finsler Geometry and Convex Geometry. As we all know, the vertical bisector is a geometric figure with simple structure. In the two-dimensional Euclidean space, the vertical bisector is just a straight line. But, in the non-Euclidean Minkowski space, the bisector is not necessarily a straight line. Therefore, based on the purpose of exploring the properties of bisector in arbitrary Banach space, this paper starts from the perpendicular bisector of Euclidean unit circle and uses the *P*-angle function to define a new angular geometric constant, as shown below:

$$S_P(X) = \sup \{\cos ang_P(x+y, x-y) : x, y \in S_X, x \neq \pm y\}.$$

Then, this paper deduces the relationship between the new geometric constant and James constant, von-Nuemann constant, the module of convexity and the module of smoothness. Finally, the judgment theorems of the geometric properties of Banach spaces are obtained, including uniform non-squareness, uniform convexity, strict convexity and uniform normal structure. In addition, we focus on the properties of Minkowski spaces by using the new constants $S_P(X)$, and derive the sufficient conditions for the fixed point property and super-normal structure.

## 2 Notations and Preliminaries

Throughout the paper, let *X* be real Banach space with *dimX* > 2. The unit ball and the unit sphere of *X* are denoted by $B_X$ and $S_X$, respectively. Now, let's recall some concepts of geometric properties of Banach spaces and significant functions that we need in this paper.

**Definition 2.1.** *[7] The Banach space X is said to be uniformly convex whenever given $0 < \varepsilon \leq 2$, there exists $\delta > 0$ such that if $x, y \in S_X$ and $\|x - y\| \geq \varepsilon$, then $\frac{\|x-y\|}{2} \leq 1 - \delta$.*

In order to explore the convexity of Banach spaces, Clarkson [11, 15] defined the module of convexity, and derived sufficient conditions for uniform convexity and strict convexity.

**Definition 2.2.** *[11, 15] Let X be Banach space, then the module of convexity is defined by*

$$\delta_X(\varepsilon) = \inf\left\{1 - \frac{\|x+y\|}{2} : x, y \in S_X, \|x-y\| = \varepsilon\right\}$$

where $\varepsilon \in [0, 2]$.

And the characteristic of convexity of $(X, \|\cdot\|)$ is defined by $\varepsilon_0(X) = \sup\{\varepsilon \in [0, 2] : \delta_X(\varepsilon) = 0\}$.

**Definition 2.3.** *[13] The Banach space X is called to be uniformly non-square if there exists $\delta \in (0, 1)$ such that for any $x, y \in S_X$, either $\frac{\|x-y\|}{2} \leq 1 - \delta$ or $\frac{\|x+y\|}{2} \leq 1 - \delta$.*

In order to explore the non squareness of Banach spaces, Gao [13] and Clarkson [6] defined the James constant and the von-Neumann constant respectively, as shown below:

**Definition 2.4.** *[13] Let X be Banach space, then the James constant is defined by*

$$J(X) = \sup\{\min\{\|x+y\|, \|x-y\|\} : x, y \in S_X\}.$$

**Definition 2.5.** *[6] Let X be Banach space, then the von-Neumann constant is defined by*

$$C_{NJ}(X) = \sup\left\{\frac{\|x+y\|^2 + \|x-y\|^2}{2\|x\|^2 + 2\|y\|^2} : x, y \in X, (x,y) \neq (0,0)\right\}.$$

*And the modified von-Neumann constant is defined by*

$$C'_{NJ}(X) = \sup\left\{\frac{\|x+y\|^2 + \|x-y\|^2}{4} : x, y \in S_X\right\}.$$

**Definition 2.6.** *[19] Let X be Banach space, then the function $\gamma_X(t) : [0, 1] \to [0, 4]$ is defined by*

$$\gamma_X(t) = \sup\left\{\frac{\|x+ty\|^2 + \|x-ty\|^2}{2} : x, y \in S_X\right\}.$$

In addition, in order to better characterize the properties of Banach spaces, Zbaganu generalized the constant $C_{NJ}(X)$ in 2001 and introduced the following constant:

**Definition 2.7.** *[8] Let X be Banach space, then the Zbaganu constant is defined by*

$$C_Z(X) = \sup\left\{\frac{\|x+y\|\|x-y\|}{\|x\|^2 + \|y\|^2} : x, y \in X, (x, y) \neq (0, 0)\right\}.$$

Alonso and Martin [18] proved the existence of Banach space $X$ such that $C_Z(X) < C_{NJ}(X)$.

**Definition 2.8.** *[7] The Banach space X is called to be strictly convex if $\|x\|=\|y\|=1$ and $x \neq y$ imply $\|x+y\|<2$.*

**Definition 2.9.** *[3, 16] Let X be Banach space and A be the bounded closed convex subset of X, then diamA $= \sup\{\|x-y\| : x, y \in A\}$ is called as the diameter of A and $r(A) = \inf\{\sup\{\|x-y\|\} : y \in A\}$ is called as the Chebyshev radius of A. The Banach space X is said to have normal structure if $r(A) <$ diamA for every bounded closed convex subset A of X with diamA $> 0$. The Banach space X is said to have uniform normal structure if*

$$\inf\left\{\frac{diam A}{r(A)}\right\} > 1$$

*with diamA $> 0$. Let $\mathcal{U}$ be a free ultrafilter over N (the set of natural integers), then the ultraproduct $\mathcal{X}$ of X is the quotient space of*

$$l_\infty(X) = \left\{(x_n); x_n \in X \text{ and } \|(x_n)\| = \sup_n \|x_n\| < \infty\right\}$$

*by*

$$\mathcal{N} = \left\{(x_n) \in l_\infty(X); \lim_{n \to \mathcal{U}} \|x_n\| = 0\right\}.$$

*The Banach space X is said to have super-normal structure if any ultraproduct X of X has normal structure.*

**Definition 2.10.** *[9] The Banach space X is called to be uniformly smooth if for any $\varepsilon > 0$, there exists $\eta > 0$ such that $1 - \frac{\|x-y\|}{2} \leq \varepsilon \|x - y\|$ for any $x, y \in S_X$ and $\|x - y\| \leq \eta$.*

In order to explore the uniform smoothness of Banach spaces, Day [9] introduced the module of smoothness and derived some characterizations of the smoothness of Banach spaces.

**Definition 2.11.** *[9] Let X be Banach space, then the module of smoothness $\rho_X(t)$ is defined by*

$$\rho_X(t) = \sup\left\{\frac{\|x+ty\| + \|x-ty\|}{2} - 1 : x, y \in S_X,\right\},$$

*where $t \in [0, +\infty)$.*

Next, we list some conclusions about the geometric properties of Banach spaces as follows:

**Lemma 2.12.** *Let X be Banach space.*
(i) *If $\varepsilon_0(X) < 2$, then X has the fixed point property [11].*
(ii) *$\sqrt{2} \leq J(X) \leq 2$ [12, 14].*
(iii) *X is uniformly non-square if and only if $J(X) < 2$ [12].*

(iv) *If $J(X) < \frac{1+\sqrt{5}}{2}$, then X has uniform normal structure [10]*.

(v) *If $2\gamma_X(t) < 1 + (1 + t)^2$ for some $t \in (0, 1]$, then X has super-normal structure [19]*.

(vi) *X is Hilbert space if and only if $\gamma_X(t) = 1 + t^2$ for any $t \in [0, 1]$ [19]*.

As we all know, for general normed spaces, the parallelogram rule can describe inner product spaces. In [4], this rule is extended to the following form:

**Lemma 2.13.** *[4] Let X be real normed linear space, then $(X, \|\cdot\|)$ is an inner product space if and only if for any $x, y \in S_X$, there exist $\alpha, \beta \neq 0$ such that $\|\alpha x + \beta y\|^2 + \|\alpha x - \beta y\|^2 \sim 2(\alpha^2 + \beta^2)$, where $\sim$ stands for =, $\leq$ or $\geq$.*

In [5], Vitor Balestro defined the *P*-angle function from Pythagorean orthogonality, which extended the original concept of angle in Euclidean spaces, as shown below:

**Definition 2.14.** *[5] Let X be Banach space, $x, y \in X_o = X - \{o\}$, then the P−angle between two non-zero vectors x, y to be*

$$ang_P(x, y) = \arccos\left(\frac{\|x\|^2 + \|y\|^2 - \|x-y\|^2}{2\|x\|\|y\|}\right).$$

## 3 The new angular geometric constant $S_P(X)$

Motivated by the form of $J(X)$, we can know that this constant is defined by using the geometric means of the variable lengths of the sides of triangles with vertices $x, -x$ and $y$, where $x, y$ are points on the unit sphere of the Banach space $X$, which can be explained from the following images:

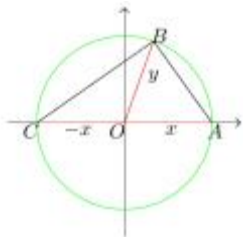

Figure 1:

As shown in the Figure 1, for the two-dimensional Euclidean plane, we have

$$\overrightarrow{BA} = x - y, \overrightarrow{BC} = x + y,$$

then the geometric constant $J(X)$ is defined by comparing the lengths of line segments $BC$ and $AB$. In 2008, by using the same ideas, Alonso and Llorens-Fuster [1] defined two geometric constants $t(X), T(X)$, as follows:

$$t(X) = \inf_{x \in S_X} \sup_{y \in S_X} \sqrt{\|x+y\|\|x-y\|}, T(X) = \sup_{x,y \in S_X} \sqrt{\|x+y\|\|x-y\|}.$$

In view of the geometric constants, there are many papers related to length comparison. In fact, in the Figure 1, according to the relationship between the center angle and the circumference angle in the elementary mathematical knowledge, it is easy to know that $\cos ang_P(x + y, x − y) = 0$. Thus, by changing the comparison of the lengths of the two line segments to explore the *P*−angle between $\overrightarrow{AB}$ and $\overrightarrow{CB}$ in the Figure 1, we obtain the following definition of the new angular geometric constant:

**Definition 3.1.** *Let X be Banach space, then the P-angle constant $S_P(X)$ is defined by*

$$S_P(X) = \sup\{\cos ang_P(x+y, x-y) : x, y \in S_X, x \neq \pm y\}.$$

**Remark 3.2.** *In fact, for two-dimensional Euclidean plane, the geometric constant $S_P(X)$ has strong geometric intuition, which can be explained from the following images:*

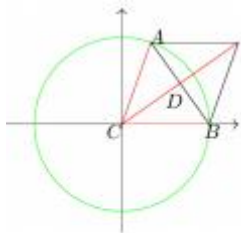

Figure 2:

As shown in the Figure 2, in the arbitrary triangle ABC of the two-dimensional Euclidean plane unit circle, we set $\overrightarrow{CA} = \beta, \overrightarrow{CB} = \alpha$, then $\overrightarrow{AB} = \alpha - \beta, \overrightarrow{CD} = \frac{1}{2}(\alpha + \beta)$. According to the properties of isosceles triangles, it can be seen that the vector $\overrightarrow{AB}$ is perpendicular to the vector $\overrightarrow{CD}$, that is, the P−angle between the two vectors $\alpha + \beta, \alpha - \beta$ is equal to $\frac{\pi}{2}$, thus $\cos ang_p(\alpha + \beta, \alpha - \beta) = 0$, which implies that $S_P(X) = 0$.

**Proposition 3.3.** *Let X be Banach space, then $0 \leq S_P(X) \leq \frac{1}{2}$.*

*Proof.* For the any normed space X, there exist $x_0, y_0 \in S_X$ such that $\|x_0+y_0\|=\|x_0-y_0\|=\sqrt{2}$ [2, Proposition 5], thus $S_P(X) \geq \cos ang_P(x_0+y_0, x_0-y_0) = 0$.

In addition, letting $x, y \in S_X$ and $x \neq \pm y$, we can assume $\|x+y\| \leq \|x-y\|$, then

$$\|x+y\|^2 + \|x-y\|^2 - 4 - \|x+y\|\|x-y\| = \left(\|x+y\| - \frac{\|x-y\|}{2}\right)^2 + \frac{3\|x-y\|^2}{4} - 4$$
$$\leq \left(\|x-y\| - \frac{\|x-y\|}{2}\right)^2 + \frac{3\|x-y\|^2}{4} - 4$$
$$= \|x-y\|^2 - 4 \leq 0,$$

which implies $S_P(X) \leq \frac{1}{2}$.

**Example 3.4.** *Let $l_p$ ( $1 \leq p < \infty$ ) be the linear space of all sequences in R such that $\sum_{i=1}^{\infty} |x_i|^p < \infty$ with the norm defined by*

$$\|x\|_p = \left(\sum_{i=1}^{\infty} |x_i|^p\right)^{\frac{1}{p}}.$$

*(i) If $2 \leq p < \infty$, we choose* $x_0 = \left(\frac{1}{2^{\frac{1}{p}}}, \frac{1}{2^{\frac{1}{p}}}, 0, \cdots\right), y_0 = \left(\frac{1}{2^{\frac{1}{p}}}, -\frac{1}{2^{\frac{1}{p}}}, 0, \cdots\right)$, then

$$\|x_0 + y_0\| = \|x_0 - y_0\| = 2^{1-\frac{1}{p}},$$

*that is,* $S_P(l_p) \geq \cos ang_P(x_0, y_0) = \frac{2^{2-\frac{2}{p}} + 2^{2-\frac{2}{p}} - 4}{2 \cdot 2^{2-\frac{2}{p}}} = 1 - 2^{\frac{2}{p}-1}$.

*(ii) If $1 \leq p < 2$, we choose $x_0 = (1, 0, \cdots), y_0 = (0, 1, 0, \cdots)$, then*

$$\|x_0 + y_0\| = \|x_0 - y_0\| = 2^{\frac{1}{p}},$$

*that is,* $S_P(l_p) \geq \cos ang_P(x_0, y_0) = \frac{2^{\frac{2}{p}} + 2^{\frac{2}{p}} - 4}{2 \cdot 2^{\frac{2}{p}}} = 1 - 2^{1-\frac{2}{p}}.$

*Therefore* $S_P(l_p) \geq 1 - 2^{-|\frac{2}{p} - 1|}.$

In addition to estimating the upper and lower bounds of the constant $S_P(X)$, we can also give the characterization of the inner product space with the new constant

**Theorem 3.5.** *Let X be Banach space, then X is Hilbert space if and only if $S_P(X) = 0$.*

*Proof.* Assume that $X$ is Hilbert space, then $\|x+y\|^2 + \|x-y\|^2 = 2\|x\|^2 + 2\|y\|^2$ for $x, y \in X$.

Letting $x, y \in S_X$ and $x \neq \pm y$, then $\cos ang_P(x+y, x-y) = 0$, thus $S_P(X) = 0$.

Assume $S_P(X) = 0$, then

$$\frac{\|x+y\|^2 + \|x-y\|^2 - 4}{2\|x+y\|\|x-y\|} \leq 0$$

for any $x, y \in S_X$ and $x \neq \pm y$, therefore $\|x+y\|^2 + \|x-y\|^2 \geq 4$ for any $x, y \in S_X$, which implies that $X$ is Hilbert space.

## 4 Some inequalities for $S_P(X)$ and significant constants

In this section, we will study some inequalities for $S_P(X)$ and some geometric constants, including the James constant $J(X)$, the modified von-Neumann constant $C_{NJ}'(X)$, the module of convexity $\delta_X(\varepsilon)$ and the module of smoothness $\rho_X(t)$. Moreover, these inequalities will help us to discuss the relations between $S_P(X)$ and some properties of Banach spaces in the next section.

**Theorem 4.1.** *Let X be Banach space, then* $S_P(X) \geq \frac{1}{2}(C_{NJ}'(X) - 1).$

*Proof.* According to the definition of $S_P(X)$, we can get that

$$S_P(X) \geq \cos ang_P(x+y, x-y) = 2 \cdot \frac{\frac{\|x+y\|^2 + \|x-y\|^2}{4} - 1}{\|x+y\|\|x-y\|}$$

for any $x, y \in S_X$ and $x \neq \pm y$, thus $\frac{\|x+y\|\|x-y\|}{2} S_P(X) \geq \frac{\|x+y\|^2 + \|x-y\|^2}{4} - 1.$

Note that $\|x+y\|\|x-y\| \leq 2\min\{\|x+y\|, \|x-y\|\}$, then

$$J(X) S_P(X) \geq \frac{\|x+y\|\|x-y\|}{2} S_P(X) \geq \frac{\|x+y\|^2 + \|x-y\|^2}{4} - 1$$

which shows $2S_P(X) \geq J(X) S_P(X) \geq C_{NJ}'(X) - 1$, thus $S_P(X) \geq \frac{1}{2}(C_{NJ}'(X) - 1).$

**Corollary 4.2.** *Let X be Banach space, then* $S_P(X) \leq \frac{C_{NJ}'(X) - 1}{\delta_X(0)}.$

*Proof.* (i) If $S_P(X) = 0$, the conclusion is clearly valid.

(ii) If $S_P(X) > 0$, then for any $\varepsilon \in (0, S_P(X))$, there exist $x_0, y_0 \in S_X$ and $x \neq \pm y_0$, such that

$$S_P(X) - \varepsilon \leq \cos ang_P(x_0 + y_0, x_0 - y_0).$$

Note that $\cos ang_P(x_0 + y_0, x_0 - y_0) = 2 \cdot \frac{\frac{\|x_0+y_0\|^2+\|x_0-y_0\|^2}{4}-1}{\|x_0+y_0\|\|x_0-y_0\|}$, then

$$\frac{\|x_0 + y_0\|\|x_0 - y_0\|}{2}(S_P(X) - \varepsilon) \leq \frac{\|x_0 + y_0\|^2 + \|x_0 - y_0\|^2}{4} - 1.$$

Assume $\|x_0+y_0\| \geq \|x_0-y_0\|$, we can get that

$$1 - \frac{\|x_0+y_0\|}{2} \geq \delta_X(\|x_0 - y_0\|),$$

then

$$\|x_0 + y_0\|\|x_0 - y_0\| = \|x_0 + y_0 - 2y_0\|\|x_0 + y_0\| \geq (2 - \|x_0 + y_0\|)\|x_0 + y_0\|$$
$$\geq 2\delta_X(\|x_0 - y_0\|)\max\{\|x_0 + y_0\|, \|x_0 - y_0\|\}$$
$$\geq 2\delta_X(\|x_0 - y_0\|)S(X).$$

Since the module of convexity is non-decreasing and $\|x_0-y_0\|>0$, then $\delta_X(\|x_0-y_0\|) \geq \delta_X(0)$, which shows $\|x_0+ y_0\|\|x_0-y_0\| \geq 2\delta_X(0)S(X)$. Thus

$$\delta_X(0)S(X)(S_P(X) - \varepsilon) \leq \frac{\|x_0 + y_0\|\|x_0 - y_0\|}{2}(S_P(X) - \varepsilon)$$
$$\leq \frac{\|x_0 + y_0\|^2 + \|x_0 - y_0\|^2}{4} - 1 \leq C'_{NJ}(X) - 1.$$

Letting $\varepsilon \to 0^+$, then $\delta_X(0)S(X)S_P(X) \leq C_{NJ}'(X) - 1$. Since $S(X)J(X) = 2$ [14, Theorem 2], then $2\delta_X(0)S_P(X) \leq (C_{NJ}'(X) - 1)J(X) \leq 2(C_{NJ}'(X) - 1)$. Thus $S_P(X) \leq \frac{C'_{NJ}(X)-1}{\delta_X(0)}$.

**Theorem 4.3.** *Let X be Banach space, then* $S_P(X) \leq \frac{4\rho_X^2(1)+4\rho_X(1)-3}{8\delta_X(0)}$.

*Proof.* (i) If $S_P(X) = 0$, the conclusion is clearly valid.

(ii) If $S_P(X) > 0$, then for any $\varepsilon \in (0, S_P(X))$, there exist $x_0, y_0 \in S_X$ and $x_0 \neq y_0$, such that

$$S_P(X) - \varepsilon \leq \cos ang_P(x_0 + y_0, x_0 - y_0).$$

Note that $\cos ang_P(x_0 + y_0, x_0 - y_0) = 2 \cdot \frac{\frac{\|x_0+y_0\|^2+\|x_0-y_0\|^2}{4}-1}{\|x_0+y_0\|\|x_0-y_0\|}$, then

$$\frac{\|x_0 + y_0\|\|x_0 - y_0\|}{2}(S_P(X) - \varepsilon) \leq \frac{\|x_0 + y_0\|^2 + \|x_0 - y_0\|^2}{4} - 1.$$

Since $\|x_0+y_0\|^2 + \|x_0-y_0\|^2 \leq (\|x_0+y_0\|+\|x_0-y_0\|)^2 \leq (2\rho_X(1)+1)^2$, which implies that

$$\frac{\|x_0 + y_0\|^2 + \|x_0 - y_0\|^2}{4} - 1 \leq \frac{1}{4}(2\rho_X(1) + 1)^2 - 1 = \rho_X^2(1) + \rho_X(1) - \frac{3}{4}.$$

Since $\|x_0 + y_0\|\|x_0 - y_0\| \geq 2\delta_X(0)S(X) = \frac{4\delta_X(0)}{J(X)}$, then

$$\frac{4\delta_X(0)}{J(X)}(S_P(X) - \varepsilon) \leq \frac{\|x_0 + y_0\|\|x_0 - y_0\|}{2}(S_P(X) - \varepsilon)$$

and

$$\frac{\|x_0 + y_0\|\|x_0 - y_0\|}{2}(S_P(X) - \varepsilon) \leq \frac{\|x_0 + y_0\|^2 + \|x_0 - y_0\|^2}{4} - 1,$$

that is, $\frac{4\delta_X(0)}{J(X)}(S_P(X) - \varepsilon) \leq \rho_X^2(1) + \rho_X(1) - \frac{3}{4}$. Letting $\varepsilon \to 0^+$, then

$$4\delta_X(0)S_P(X) \le \left[\rho_X(1)^2 + \rho_X(1) - \tfrac{3}{4}\right] J(X) \le 2\rho_X^2(1) + 2\rho_X(1) - \tfrac{3}{2}.$$

This completes the proof.

**Theorem 4.4.** *Let X be Banach space, then* $S_P(X) \le \frac{\gamma_X(t) - 2t^2}{2t^2 \delta_X(0)}$, *where* $t \in (0, 1]$.

*Proof.* Letting $x, y \in S_X$ and $x \ne \pm y$, since

$$\|x+y\| = \left\|\frac{1+t}{2t} \cdot (x+ty) + \frac{t-1}{2t} \cdot (x-ty)\right\| \le \frac{1+t}{2t}\|x+ty\| + \frac{1-t}{2t}\|x-ty\|$$

and

$$\|x-y\| = \left\|\frac{1+t}{2t} \cdot (x-ty) + \frac{t-1}{2t} \cdot (x+ty)\right\| \le \frac{1+t}{2t}\|x-ty\| + \frac{1-t}{2t}\|x+ty\|,$$

then

$$2\|x+y\|\|x-y\| \cos ang_P(x+y, x-y) = \|x+y\|^2 + \|x-y\|^2 - 4$$

$$\le \left(\frac{1+t^2}{2t^2} + \frac{1-t^2}{2t^2}\right)(\|x+ty\|^2 + \|x-ty\|^2) - 4$$

$$= \frac{1}{t^2}(\|x+ty\|^2 + \|x-ty\|^2) - 4$$

$$\le \frac{2\gamma_X(t)}{t^2} - 4.$$

Since $\|x+y\|\|x-y\| \ge \frac{4\delta_X(0)}{J(X)}$, then $\frac{8\delta_X(0)}{J(X)} \cos ang_P(x+y, x-y) \le \frac{2\gamma_X(t)}{t^2} - 4$, which implies that $2\delta_X(0)S_P(X) \le \left(\frac{\gamma_X(t)}{2t^2} - 1\right) J(X) \le \frac{\gamma_X(t)}{t^2} - 2$. Thus $S_P(X) \le \frac{\gamma_X(t) - 2t^2}{2t^2 \delta_X(0)}$.

**Theorem 4.5.** *Let X be Banach space, then* $S_P(X) \le \frac{C_Z(X)}{2\delta_X(0)}$.

*Proof.* Since $\delta_X(0) \ne 0$, then let $x, y \in S_X$ and $x \ne y$. Assume $\|x-y\| \le \|x+y\|$, then

$$\|x+y\|^2 + \|x-y\|^2 - 4 \le \|x-y\|^2 \le \|x+y\|\|x-y\|.$$

Since $\|x+y\|\|x-y\| \ge 2\delta_X(0)S(X) = \frac{4\delta_X(0)}{J(X)}$, then

$$\frac{\|x+y\|^2 + \|x-y\|^2 - 4}{2\|x+y\|\|x-y\|} \le \frac{\|x+y\|\|x-y\|}{\frac{8\delta_X(0)}{J(X)}} = \frac{J(X)}{4\delta_X(0)} \cdot \frac{\|x+y\|\|x-y\|}{\|x\|^2 + \|y\|^2} \le \frac{J(X)C_Z(X)}{4\delta_X(0)},$$

which implies that $4\delta_X(0)S_P(X) \le J(X)C_Z(X) \le 2C_Z(X)$. Thus $S_P(X) \le \frac{C_Z(X)}{2\delta_X(0)}$.

Next, we focus on the relationship between $S_P(X)$ of the finite dimensional Banach spaces $X$ and other geometric constants, and then pave the way for the next part to explain the advantages of the constant $S_P(X)$ in describing the geometric properties of Banach spaces.

**Corollary 4.6.** *Let X be finite dimensional Banach space, then* $S_P(X) \ge \frac{\gamma_X(t) + t^2 - 3}{2 + 2t^2}$, *where* $t \in [0, 1]$.

*Proof.* Note that the conclusion is obviously valid for $t = 0$, so we only need to consider the case of $t \in (0, 1]$.

(i) If $X$ is Hilbert space, then $S_P(X) = 0$ and $\gamma_X(t) = 1 + t^2$ for $t \in [0, 1]$, which shows that

$$\frac{\gamma_X(t) + t^2 - 3}{2 + 2t^2} = \frac{t^2 - 1}{t^2 + 1} \leq 0 = S_P(X),$$

then the conclusion is established.

(ii) If $X$ is not Hilbert space, then we can deduce $\gamma_X(t) \neq 1 + t^2$ for any $t \in (0, 1]$ by Lemma 2.12.

In addition, there exist $x_{n,t}, y_{n,t} \in S_X$ such that

$$\lim_{n \to \infty} \frac{\|x_{n,t} + ty_{n,t}\|^2 + \|x_{n,t} - ty_{n,t}\|^2}{2} = \gamma_X(t).$$

Since $X$ is finite dimensional and $x_{n,t}, y_{n,t} \in S_X$, then there exist $x_{0,t}, y_{0,t} \in S_X$ such that

$$\lim_{n \to \infty} \|x_{n,t}\| = \|x_{0,t}\|, \lim_{n \to \infty} \|y_{n,t}\| = \|y_{0,t}\|,$$

that is, $\|x_{0,t} + ty_{0,t}\|^2 + \|x_{0,t} - ty_{0,t}\|^2 = 2\gamma_X(t)$.

If $x_{0,t} = y_{0,t}$ or $x_{0,t} = -y_{0,t}$, then $\gamma_X(t) = \frac{(1+t)^2 + (1-t)^2}{2} = 1 + t^2$, this is contradictory, thus $x_{0,t} \neq \pm y_{0,t}$.

Note that

$$\|x_{0,t} + ty_{0,t}\| = \left\|\frac{1+t}{2}(x_{0,t} + y_{0,t}) + \frac{1-t}{2}(x_{0,t} - y_{0,t})\right\| \leq \frac{1+t}{2}\|x_{0,t} + y_{0,t}\| + \frac{1-t}{2}\|x_{0,t} - y_{0,t}\|$$

and

$$\|x_{0,t} - ty_{0,t}\| = \left\|\frac{1+t}{2}(x_{0,t} - y_{0,t}) + \frac{1-t}{2}(x_{0,t} + y_{0,t})\right\| \leq \frac{1+t}{2}\|x_{0,t} - y_{0,t}\| + \frac{1-t}{2}\|x_{0,t} + y_{0,t}\|,$$

then

$$\|x_{0,t} + ty_{0,t}\|^2 + \|x_{0,t} - ty_{0,t}\|^2$$
$$\leq \frac{1+t^2}{2}(\|x_{0,t} + y_{0,t}\|^2 + \|x_{0,t} - y_{0,t}\|^2) + (1 - t^2)\|x_{0,t} + y_{0,t}\|\|x_{0,t} - y_{0,t}\|.$$

Since

$$\frac{1+t^2}{2}(\|x_{0,t} + y_{0,t}\|^2 + \|x_{0,t} - y_{0,t}\|^2)$$
$$= (1 + t^2)\|x_{0,t} + y_{0,t}\|\|x_{0,t} - y_{0,t}\| \cdot \frac{\|x_{0,t} + y_{0,t}\|^2 + \|x_{0,t} - y_{0,t}\|^2 - 4}{2\|x_{0,t} + y_{0,t}\|\|x_{0,t} - y_{0,t}\|} + 2 + 2t^2$$
$$\leq 2(1 + t^2)J(X)S_P(X) + 2 + 2t^2,$$

then

$$\|x_{0,t} + ty_{0,t}\|^2 + \|x_{0,t} - ty_{0,t}\|^2 \leq 2(1+t^2)J(X)S_P(X) + 2(1-t^2)J(X) + 2 + 2t^2,$$

which implies that $\gamma_X(t) \leq (1 + t^2)J(X)S_P(X) + (1 - t^2)J(X) + 1 + t^2$.

Thus $S_P(X) \geq \frac{\gamma_X(t) - (1-t^2)J(X) - 1 - t^2}{(1+t^2)J(X)}$. Let $x \in S_X$, then $\gamma_X(t) \geq \frac{\|x+tx\|^2 + \|x-tx\|^2}{2} = 1 + t^2$, thus

$$\frac{\gamma_X(t) - (1-t^2)J(X) - 1 - t^2}{(1+t^2)J(X)} = \frac{\gamma_X(t) - t^2 - 1}{(1+t^2)J(X)} + \frac{t^2 - 1}{1+t^2} \geq \frac{\gamma_X(t) - t^2 - 1}{2 + 2t^2} + \frac{t^2 - 1}{1+t^2}$$
$$= \frac{\gamma_X(t) + t^2 - 3}{2 + 2t^2}.$$

Then $S_P(X) \geq \frac{\gamma_X(t) + t^2 - 3}{2 + 2t^2}$.

In order to prove the following corollary, we first give the characterization of the inner product space by using the constant $C_Z(X)$.

**Lemma 4.7.** *Let X be Banach space, then X is Hilbert space if and only if $C_Z(X) = 1$.*

*Proof.* Assume that $X$ is Hilbert space, then $\|x+y\|^2 + \|x-y\|^2 = 2(\|x\|^2 + \|y\|^2)$ for any $x, y \in X$, thus

$$2\|x+y\|\|x-y\| \leq \|x+y\|^2 + \|x-y\|^2 = 2(\|x\|^2 + \|y\|^2),$$

that is, $\|x+y\|\|x-y\| \leq \|x\|^2 + \|y\|^2$. Hence $C_Z(X) \leq 1$, which implies that $C_Z(X) = 1$.

Assume $C_Z(X) = 1$, then $\frac{\|x+y\|\|x-y\|}{\|x\|^2+\|y\|^2} \leq 1$ for $x, y \in X$ and $(x, y) \neq (0, 0)$, that is, $\|x\|^2 + \|y\|^2 \geq 4\|x+y\|\|x-y\|$ for any $x, y \in X$. Letting $u = x+y$, $v = x-y$, then $x = \frac{u+v}{2}, y = \frac{u-v}{2}$. Hence we can get

$$\|u+v\|^2 + \|u-v\|^2 \geq 4\|u\|\|v\|.$$

Therefore, we have $\|u+v\|^2 + \|u-v\|^2 \geq 4$ for any $u, v \in S_X$, that is, $X$ is Hilbert space.

**Corollary 4.8.** *Let X be finite dimensional Banach space, then* $S_P(X) \geq \frac{3\sqrt{C_Z(X)-2}}{C_Z(X)} - 1$.

*Proof.* (i) If $X$ is Hilbert space, then $S_P(X) = 0$ and $C_Z(X) = 1$, thus

$$\frac{3\sqrt{C_Z(X)-2}}{C_Z(X)} - 1 \leq S_P(X),$$

then the conclusion is clearly established.

(ii) If $X$ is not Hilbert space, then $C_Z(X) > 1$. Since $X$ is finite dimensional, then there exist $u_0, v_0 \in X$ and $(u_0, v_0) \neq (0, 0)$ such that

$$\frac{\|u_0+v_0\|\|u_0-v_0\|}{\|u_0\|^2 + \|v_0\|^2} = C_Z(X).$$

In fact, $\|u_0\|, \|v_0\| \neq 0$, since if $\|u_0\| = 0$ or $\|v_0\| = 0$, then $C_Z(X) = 1$, which is contradictory to $C_Z(X) > 1$.

Suppose $\|u_0\| \leq \|v_0\|$, and let $x_0 = \frac{u_0}{\|u_0\|}, y_0 = \frac{v_0}{\|v_0\|} \in S_X$. If $x_0 = y_0$ or $x_0 = -y_0$, then

$$u_0 = \lambda_0 v_0 \text{ or } u_0 = -\lambda_0 v_0, \text{ where } \lambda_0 = \frac{\|u_0\|}{\|v_0\|} \in (0, 1],$$

thus

$$C_Z(X) = \frac{\|u_0+v_0\|\|u_0-v_0\|}{\|u_0\|^2+\|v_0\|^2} = \frac{(\lambda_0+1)(1-\lambda_0)}{\lambda_0^2+1} = \frac{1-\lambda_0^2}{1+\lambda_0^2} < 1,$$

which is contradictory to $C_Z(X) > 1$. Hence $x_0 \neq \pm y_0$

Note that

$$\|x_0+y_0\| = \frac{1}{\|u_0\|\|v_0\|} \|\|v_0\|u_0 + \|u_0\|v_0\| = \left\|\frac{\|u_0\|+\|v_0\|}{2\|u_0\|\|v_0\|}(u_0+v_0) + \frac{\|v_0\|-\|u_0\|}{2\|u_0\|\|v_0\|}(u_0-v_0)\right\|$$

$$\leq \frac{\|u_0\|+\|v_0\|}{2\|u_0\|\|v_0\|}\|u_0+v_0\| + \frac{\|v_0\|-\|u_0\|}{2\|u_0\|\|v_0\|}\|u_0-v_0\|$$

and

$$\|x_0-y_0\| = \frac{1}{\|u_0\|\|v_0\|} \|\|v_0\|u_0 - \|u_0\|v_0\| = \left\|\frac{\|u_0\|+\|v_0\|}{2\|u_0\|\|v_0\|}(u_0-v_0) + \frac{\|v_0\|-\|u_0\|}{2\|u_0\|\|v_0\|}(u_0+v_0)\right\|$$

$$\leq \frac{\|u_0\|+\|v_0\|}{2\|u_0\|\|v_0\|}\|u_0-v_0\| + \frac{\|v_0\|-\|u_0\|}{2\|u_0\|\|v_0\|}\|u_0+v_0\|,$$

then

$$\|x_0 + y_0\|\|x_0 - y_0\| \leq \frac{\|v_0\|^2 - \|u_0\|^2}{4\|u_0\|^2\|v_0\|^2}(\|u_0+v_0\|^2 + \|u_0-v_0\|^2) + \frac{\|u_0\|^2 + \|v_0\|^2}{2\|u_0\|^2\|v_0\|^2}\|u_0+v_0\|\|u_0-v_0\|.$$

In addition, since

$$\|x_0 - y_0\| \geq \left| \frac{\|u_0\| + \|v_0\|}{2\|u_0\|\|v_0\|}\|u_0+v_0\| - \frac{\|v_0\| - \|u_0\|}{2\|u_0\|\|v_0\|}\|u_0-v_0\| \right|$$

and

$$\|x_0 - y_0\| \geq \left| \frac{\|u_0\| + \|v_0\|}{2\|u_0\|\|v_0\|}\|u_0-v_0\| - \frac{\|v_0\| - \|u_0\|}{2\|u_0\|\|v_0\|}\|u_0+v_0\| \right|,$$

then

$$\|x_0+y_0\|^2 + \|x_0-y_0\|^2 \geq \frac{\|u_0\|^2 + \|v_0\|^2}{2\|u_0\|^2\|v_0\|^2}(\|u_0+v_0\|^2 + \|u_0-v_0\|^2) - \frac{\|v_0\|^2 - \|u_0\|^2}{\|u_0\|^2\|v_0\|^2}\|u_0+v_0\|\|u_0-v_0\|.$$

Setting $A = \frac{\|v_0\|^2 - \|u_0\|^2}{4\|u_0\|^2\|v_0\|^2}$, $B = \frac{\|v_0\|^2 + \|u_0\|^2}{4\|u_0\|^2\|v_0\|^2}$, then $2A + B = \frac{1}{\|u_0\|^2}$, $B - 2A = \frac{1}{\|v_0\|^2}$ and

$$\begin{aligned}
S_P(X) &\geq \frac{\|x_0+y_0\|^2 + \|x_0-y_0\|^2 - 4}{2\|x_0+y_0\|\|x_0-y_0\|} \\
&\geq \frac{1}{2} \cdot \frac{B(\|u_0+v_0\|^2 + \|u_0-v_0\|^2) - 4A\|u_0+v_0\|\|u_0-v_0\| - 4}{A(\|u_0+v_0\|^2 + \|u_0-v_0\|^2) + B\|u_0+v_0\|\|u_0-v_0\|} \\
&= -\frac{2A}{B} + \frac{\left(\frac{2A^2}{B} + \frac{B}{2}\right)(\|u_0+v_0\|^2 + \|u_0-v_0\|^2) - 2}{A(\|u_0+v_0\|^2 + \|u_0-v_0\|^2) + B\|u_0+v_0\|\|u_0-v_0\|}.
\end{aligned}$$

Since

$$\frac{\|u_0+v_0\|^2 + \|u_0-v_0\|^2}{A(\|u_0+v_0\|^2 + \|u_0-v_0\|^2) + B\|u_0+v_0\|\|u_0-v_0\|} = \frac{1}{A + \frac{B\|u_0+v_0\|\|u_0-v_0\|}{\|u_0+v_0\|^2 + \|u_0-v_0\|^2}}$$

$$\geq \frac{1}{A + \frac{B}{2}} = \frac{2}{2A+B}$$

and

$$A(\|u_0+v_0\|^2 + \|u_0-v_0\|^2) + B\|u_0+v_0\|\|u_0-v_0\| \geq (2A+B)\|u_0+v_0\|\|u_0-v_0\|$$

$$= (2A+B)C_Z(X)\left(\frac{1}{2A+B} + \frac{1}{B-2A}\right)$$

$$= \frac{2BC_Z(X)}{B-2A},$$

then

$$S_P(X) \geq -\frac{2A}{B} + \left(\frac{2A^2}{B} + \frac{B}{2}\right) \cdot \frac{2}{2A+B} - \frac{2}{\frac{2BC_Z(X)}{B-2A}} = \frac{1-t_0}{1+t_0} + \frac{t_0 - 1}{C_Z(X)},$$

where $t_0 = \frac{2A}{B} \in (0, 1)$. Setting $f(t) = \frac{1-t}{1+t} + \frac{t-1}{C_Z(X)}, t \in [0, 1]$, by $f'(t_1) = 0$, we can get

$$t_1 = \sqrt{C_Z(X)} - 1 \in (0, 1),$$

then $\inf\{f(t) : t \in [0, 1]\} = f(t_1)$, which implies that

$$\frac{1-t_0}{1+t_0} + \frac{t_0-1}{C_Z(X)} = f(t_0) \geq f(t_1) = \frac{3\sqrt{C_Z(X)} - 2}{C_Z(X)} - 1,$$

that is, $S_P(X) \geq \frac{3\sqrt{C_Z(X)}-2}{C_Z(X)} - 1$.

# 5 The relations between $S_P(X)$ and geometric properties of Banach spaces

In this section, we will discuss the relations between $S_P(X)$ and geometric properties of Banach spaces and obtain the necessary and sufficient condition for uniform non-squareness, and the sufficient conditions for uniform convexity and the normal structure. In particular, the sufficient condition for the fixed point property is derived in the Minkowski spaces.

**Theorem 5.1.** *Let X be Banach space, then $S_P(X) < \frac{1}{2}$ if and only if X is uniformly non-square.*

*Proof.* Suppose $S_P(X) < \frac{1}{2}$, then since $C'_{NJ}(X) \geq \frac{1}{2} J(X)^2$ [1, Proposition 5], we can get

$$\frac{1}{2} > S_P(X) \geq \frac{C'_{NJ}(X) - 1}{J(X)} \geq \frac{\frac{1}{2} J(X)^2 - 1}{J(X)},$$

which shows that $\frac{1}{2} J(X)^2 - \frac{1}{2} J(X) - 1 < 0$. Hence $J(X) < 2$, that is, $X$ is uniformly non-square.

Suppose $X$ is uniformly non-square, then there exists $\delta \in (0, 1)$ such that

$$\frac{\|x+y\|}{2} \leq 1-\delta \text{ or } \frac{\|x-y\|}{2} \leq 1-\delta$$

for any $x, y \in S_X$. Then letting $x, y \in S_X$ and $x \neq \pm y$, we can deduce

$$\cos \mathrm{ang}_P(x+y, x-y) \leq \frac{\min\{\|x+y\|, \|x-y\|\}^2 + (2-2\delta)^2 - 4}{2\|x+y\|\|x-y\|}$$

$$= \frac{\min\{\|x+y\|, \|x-y\|\}^2}{2\|x+y\|\|x-y\|} + \frac{2(1-\delta)^2 - 2}{\|x+y\|\|x-y\|}$$

$$\leq \frac{1}{2} + \frac{2(1-\delta)^2 - 2}{\|x+y\|\|x-y\|},$$

that is,

$$S_P(X) \leq \frac{1}{2} + \sup\left\{\frac{2(1-\delta)^2 - 2}{\|x+y\|\|x-y\|} : x, y \in S_X, x \neq \pm y\right\} \leq \frac{1}{2}.$$

If $S_P(X) = \frac{1}{2}$, then

$$\sup\left\{\frac{2(1-\delta)^2 - 2}{\|x+y\|\|x-y\|} : x, y \in S_X, x \neq \pm y\right\} = 0,$$

thus there exist $x_n, y_n \in S_X$ and $x_n \neq \pm y_n$, such that $\lim_{n\to\infty} \frac{1}{\|x_n+y_n\|\|x_n-y_n\|} = 0$, that is,

$$\lim_{n\to\infty} \|x_n + y_n\|\|x_n - y_n\| = +\infty.$$

But $\|x_n + y_n\| \leq \|x_n\| + \|y_n\| = 2$, this is impossible. Thus $S_P(X) < \frac{1}{2}$.

**Example 5.2.** *Let $C[a, b]$ be the linear space of all real valued continuous functions on $[a, b]$ with the norm defined by $\|x\| = \sup\{|x(t)| : t \in [a, b]\}$.*

*Letting $x_0 = 1, y_0 = \frac{2}{b-a}(t-b) - 1$, then $\|x_0+y_0\| = \|x_0-y_0\| = 2$, that is*

$$S_P(C[a,b]) \geq \frac{\|x_0+y_0\|^2 + \|x_0-y_0\|^2 - 4}{2\|x_0+y_0\|\|x_0-y_0\|} = \frac{1}{2},$$

Zhijian Yang and Yongjin Li

which shows that $S_P(C[a,b]) = \frac{1}{2}$, Then $C[a, b]$ is not uniformly non-square.

**Example 5.3.** *By Example 3.4, we can deduce that $S_P(l_1) = \frac{1}{2}$, which shows $l^1$ is not uniformly non-square. In addition, let $l_\infty$ be the linear space of all bounded sequences in R with the norm defined by*

$$\|x\|_\infty = \sup\{|x_n| : 1 \leq n \leq \infty\}.$$

*We choose $x_0 = (1, 1, 0, \cdots)$, $y_0 = (1, -1, 0, \cdots)$, then $\|x_0+y_0\| = \|x_0-y_0\| = 2$, thus*

$$S_P(l_\infty) \geq \frac{\|x_0 + y_0\|^2 + \|x_0 - y_0\|^2 - 4}{2\|x_0 + y_0\|\|x_0 - y_0\|} = \frac{1}{2},$$

*that is, $S_P(l_\infty) = \frac{1}{2}$. Hence $l_\infty$ also is not uniformly non-square.*

**Theorem 5.4.** *Let X be Banach space, if $S_P(X) < \frac{3-\sqrt{5}}{4}$, then X has uniform normal structure.*

*Proof.* Since $S_P(X) < \frac{3-\sqrt{5}}{4}$, and $C'_{NJ}(X) \geq \frac{1}{2}J(X)^2$ [1, Proposition 5], then

$$\frac{3-\sqrt{5}}{4} > S_P(X) \geq \frac{C'_{NJ}(X) - 1}{J(X)} \geq \frac{\frac{1}{2}J(X)^2 - 1}{J(X)},$$

that is, $\frac{1}{2}J(X)^2 - \frac{3-\sqrt{5}}{4}J(X) - 1 < 0$. Thus

$$\left(J(X) - \frac{1+\sqrt{5}}{2}\right)\left(J(X) + \sqrt{5} - 1\right) < 0,$$

which implies $J(X) < \frac{1+\sqrt{5}}{2}$. Hence by Lemma 2.11, we know that $X$ has uniform normal structure.

The Normal structure and the fixed point property are closely related. The next corollary shows that the geometric constant $S_P(X)$ can describe the fixed point property and super-normal structure of real finite dimensional Banach spaces.

**Corollary 5.5.** *Let X be finite dimensional Banach space.*

(i) *If $S_P(X) < \frac{1}{2}$, then X has the fixed point property.*

(ii) *If $S_P(X) < \frac{1}{8}$, then X has super-normal structure.*

*Proof.* (i) Assume that $X$ doesn't have the fixed point property, then by Lemma 2.11, we can get $\varepsilon_0(X) = 2$.

Thus there exists $\varepsilon_n$ satisfying $\delta_X(\varepsilon_n) = 0$, such that $\varepsilon_n \to 2$. Since $\delta_X(\varepsilon_n) = 0$, then there exist $x_{m,n}, y_{m,n} \in S_X$ satisfying $\|x_{m,n}-y_{m,n}\| = \varepsilon_n$ such that $\lim_{m\to\infty}\left(1 - \frac{\|x_{m,n}+y_{m,n}\|}{2}\right) = 0$ which implies that $\lim_{m\to\infty}\|x_{m,n} + y_{m,n}\| = 2$.

Note that $X$ is finite dimensional, then there exist $x_{0,n}, y_{0,n} \in S_X$ such that

$$\lim_{m\to\infty}\|x_{m,n}\| = \|x_{0,n}\|, \lim_{m\to\infty}\|y_{m,n}\| = \|y_{0,n}\|.$$

Hence $\|x_{0,n} + y_{0,n}\| = 2$ and $\|x_{0,n}-y_{0,n}\| = \varepsilon_n$. Then we can deduce

$$S_P(X) \geq \cos ang_P(x_{0,n} + y_{0,n}, x_{0,n} - y_{0,n}) = \frac{4 + \varepsilon_n^2 - 4}{4\varepsilon_n} = \frac{\varepsilon_n}{4}.$$

Letting $n \to \infty$, then $S_P(X) \geq \frac{1}{2}$, which is contradictory to $S_P(X) < \frac{1}{2}$. Therefore $X$ has the fixed point property.

(ii) Since $S_P(X) \geq \frac{\gamma_X(t)+t^2-3}{2+2t^2}$, then $\frac{\gamma_X(t)+t^2-3}{2+2t^2} < \frac{1}{8}$, that is,

$$2\gamma_X(t) < -\frac{3}{2}t^2 + \frac{13}{2}.$$

Let $t_0 = 1$, since $1+(1+t_0)^2 = 5 = -\frac{3}{2}t_0^2 + \frac{13}{2}$, then $2\gamma_X(t) < 1+(1+t)^2$ for some $t \in (0, 1]$.

Hence $X$ has super-normal structure.

**Proposition 5.6.** *Let $X$ be Banach space.*

(i) *If $S_P(X) = \frac{1}{2}$, then $X$ is not uniformly convex.*

(ii) *If $S_P(X) = \frac{1}{2}$ and $X$ is finite dimensional, then $X$ is not strictly convex.*

*Proof.* (i) Since $S_P(X) = \frac{1}{2}$, then $X$ is not uniformly non-square. Thus there exist $x_0, y_0 \in S_X$ such that

$$\min\{\|x_0+y_0\|, \|x_0-y_0\|\} > 2-2\delta$$

for any $\delta \in (0, 1]$. Assume $\|x_0+y_0\| \geq \|x_0-y_0\|$, then $\|x_0+y_0\| \geq \|x_0-y_0\| > 2-2\delta$.

Letting $\varepsilon_0 = 2-2\delta \in (0,2]$, thus there exist $x_0, y_0 \in S_X$ satisfying $\|x_0-y_0\| > \varepsilon_0$ such that $\|x_0+y_0\| > 2-2\delta_0$ for any $\delta_0 = \delta \in (0,1]$. Hence $X$ is not uniformly convex.

(ii) Since $S_P(X) = \frac{1}{2}$ and $X$ is finite dimensional, then there exist $x_0, y_0 \in S_X$, such that

$$\frac{\|x_0+y_0\|^2 + \|x_0-y_0\|^2 - 4}{2\|x_0+y_0\|\|x_0-y_0\|} = \frac{1}{2},$$

that is, $4 = \|x_0+y_0\|^2 + \|x_0-y_0\|^2 - \|x_0+y_0\|\|x_0-y_0\|$.

Assume that $\|x_0+y_0\| \geq \|x_0-y_0\|$, then

$$4 = \|x_0+y_0\|^2 + \|x_0-y_0\|^2 - \|x_0+y_0\|\|x_0-y_0\|$$
$$= \left(\|x_0-y_0\| - \frac{\|x_0+y_0\|}{2}\right)^2 + \frac{3\|x_0+y_0\|^2}{4}$$
$$\leq \left(\|x_0+y_0\| - \frac{\|x_0+y_0\|}{2}\right)^2 + \frac{3\|x_0+y_0\|^2}{4} = \|x_0+y_0\|^2 \leq 4,$$

which implies that $\|x_0+y_0\| = \|x_0-y_0\| = 2$, thus $X$ is not strictly convex.

**Data Availability**

No data were used to support this study.

**Conflicts of Interest**

The author(s) declare(s) that there is no conflict of interest regarding the publication of this paper.

**Funding Statement**

This work was supported by the National Natural Science Foundation of P. R. China (Nos. 11971493 and 12071491).